\documentclass[a4paper,12pt]{amsart}
\usepackage{amssymb, mathtools}
\usepackage{graphicx}
\usepackage{ifthen}
\nonstopmode \numberwithin{equation}{section}
\newtheorem{thm}{Theorem}
\newtheorem{cor}{Corollary}
\newtheorem{lem}{Lemma}

\newtheorem{conj}{Conjecture}

\theoremstyle{definition}
\newtheorem{defn}{Definition}[section]

\newtheorem{prob}[equation]{Problem}

\newenvironment{rem}{%
\bigskip
\noindent \textsl{{\sl Remark. }}}{\bigskip}
\newenvironment{rems}{%
\bigskip
\noindent \textsl{{\sl Remarks. }}}{\bigskip}

\newcounter {own}
\def\theown {\thesection       .\arabic{own}}

\newenvironment{pf}[1][]{%
 \vskip 3mm
 \noindent
 \ifthenelse{\equal{#1}{}}%
  {{\slshape Proof. }}%
  {{\slshape #1.} }%
 }%
{\qed\bigskip}

\newcounter{alphabet}
\newcounter{tmp}
\newenvironment{Thm}[1][]{\refstepcounter{alphabet}%
\bigskip%
\noindent%
{\bf Theorem \Alph{alphabet}}%
\ifthenelse{\equal{#1}{}}{}{ (#1)}%
{\bf .} \itshape}{\vskip 8pt}

\newcommand{\B}{{\mathcal B}}

\newcommand{\ID}{{\mathbb D}}
\newcommand{\IN}{{\mathbb N}}
\newcommand{\IC}{{\mathbb C}}
\newcommand{\K}{{\mathcal K}}

\newcommand{\IB}{{\mathcal B}}
\newcommand{\IH}{{\mathcal H}}

\def\be{\begin{equation}}
\def\ee{\end{equation}}

\newcommand{\bee}{\begin{enumerate}}
\newcommand{\eee}{\end{enumerate}}

\newcommand{\blem}{\begin{lem}}
\newcommand{\elem}{\end{lem}}
\newcommand{\bthm}{\begin{thm}}
\newcommand{\ethm}{\end{thm}}
\newcommand{\bcor}{\begin{cor}}
\newcommand{\ecor}{\end{cor}}
\newcommand{\beg}{\begin{examp}}
\newcommand{\eeg}{\end{examp}}
\newcommand{\begs}{\begin{examples}}
\newcommand{\eegs}{\end{examples}}
\newcommand{\bdefe}{\begin{defn}}
\newcommand{\edefe}{\end{defn}}
\newcommand{\bprob}{\begin{prob}}
\newcommand{\eprob}{\end{prob}}
\newcommand{\bei}{\begin{itemize}}
\newcommand{\eei}{\end{itemize}}

\newcommand{\bcon}{\begin{conj}}
\newcommand{\econ}{\end{conj}}
\newcommand{\bcons}{\begin{conjs}}
\newcommand{\econs}{\end{conjs}}
\newcommand{\bprop}{\begin{propo}}
\newcommand{\eprop}{\end{propo}}
\newcommand{\br}{\begin{rem}}
\newcommand{\er}{\end{rem}}
\newcommand{\brs}{\begin{rems}}
\newcommand{\ers}{\end{rems}}
\newcommand{\bo}{\begin{obser}}
\newcommand{\eo}{\end{obser}}
\newcommand{\bos}{\begin{obsers}}
\newcommand{\eos}{\end{obsers}}
\newcommand{\bpf}{\begin{pf}}
\newcommand{\epf}{\end{pf}}
\newcommand{\ba}{\begin{array}}
\newcommand{\ea}{\end{array}}
\newcommand{\beq}{\begin{eqnarray}}
\newcommand{\beqq}{\begin{eqnarray*}}
\newcommand{\eeq}{\end{eqnarray}}
\newcommand{\eeqq}{\end{eqnarray*}}

\newcommand{\ov}{\overline}

\newcounter{minutes}
\divide\time by 60
\newcounter{hours}
\multiply\time by 60 \addtocounter{minutes}{-\time}
\begin{document}
\title{Bohr phenomenon for operator valued functions with fixed initial coefficient}
\begin{center}
{\tiny \texttt{FILE:~\jobname .tex,
        printed: \number\year-\number\month-\number\day,
        \thehours.\ifnum\theminutes<10{0}\fi\theminutes}
}
\end{center}
\author{Bappaditya Bhowmik${}^{\mathbf{*}}$}
\address{Bappaditya Bhowmik, Department of Mathematics,
Indian Institute of Technology Kharagpur, Kharagpur - 721302, India.}
\email{bappaditya@maths.iitkgp.ac.in}
\author{Nilanjan Das}
\address{Nilanjan Das, Department of Mathematics,
Indian Institute of Technology Kharagpur, Kharagpur - 721302, India.}
\email{nilanjan@iitkgp.ac.in}

\keywords{Bohr radius, Operator valued functions, Absolute values of operators.\newline
${}^{\mathbf{*}}$ Corresponding author}

\begin{abstract}
The purpose of this article is to study Bohr inequalities involving the absolute values of the
coefficients of an operator valued function. To be more specific, we establish an operator valued
analogue of a classical result regarding the Bohr phenomenon for scalar valued functions with fixed initial
coefficient. Apart from that, operator valued versions of other related and well known results are obtained.
\end{abstract}
\thanks{The first author of this article would like to
thank SERB, DST, India (Ref.No.-MTR/2018/001176) for its financial support through MATRICS grant.}

\maketitle
\pagestyle{myheadings}
\markboth{B. Bhowmik, N. Das}{Bohr phenomenon for operator valued functions with fixed initial coefficient}

\bigskip
\section{Introduction and main results}
In 1914, Harald Bohr proved the following remarkable result \cite{Bohr}.

\begin{Thm}\label{TheoA}
Let $f(z)=\sum_{n=0}^\infty a_nz^n$ be a holomorphic self mapping of the open
unit disk $\ID$. Then
\be\label{P5eq17}
\sum_{n=0}^\infty|a_n|r^n\leq 1
\ee
for all $z\in\ID$ with $|z|=r\leq 1/6$.
\end{Thm}
Inequalities of the type $(\ref{P5eq17})$ are commonly known as
\emph{Bohr inequalities} nowadays, and appearance of such inequalities in a result
is generally termed as the \emph{Bohr phenomenon}.
The above theorem was re-established with the best possible constant $1/3$ instead of $1/6$ by
Wiener, Riesz and Schur independently. Wiener's own argument can be found in \cite{Bohr} itself,
and for alternative proofs of the above inequality we refer the articles \cite{Paul1, Sid, Tom}.
Initially Theorem~A was viewed only as an outcome of the research on the absolute convergence
problem for Dirichlet series, and did not call much attention. However, in the last 25 years it
has turned to a thriving area of investigation. In fact, interest in Bohr radius problem was revived
after it found an application to the characterization problem of Banach algebras satisfying von Neumann
inequality (cf. \cite{Dix}). Since then, Bohr phenomenon is extended to multidimensional framework
(see f.i. \cite{Aiz1, Bay, Boas, Pop}), to abstract settings (cf. \cite{Aiz, Ayt, Ham}), to ordinary and
vector valued Dirichlet series (see \cite{Bala, Def1, Def2}), to free holomorphic functions and on polydomains
(cf. \cite{Pop1, Pop2}), and in numerous other directions. The readers are urged to look at the references of
the aforementioned papers to get a more complete picture of the current development in this area.

We now present two different approaches of generalizing the classical Bohr inequality, i.e. inequality (\ref{P5eq17}), which we
aim to blend together in this article. The first one was considered by Enrico Bombieri (cf. \cite[Teorema A]{Bo} or
the exposition from \cite[Chapter 8, Cor. 8.6.8]{Gar}), who established the following result.

\begin{Thm}\label{TheoB}
Let $f(z)=\lambda+\sum_{n=1}^\infty a_nz^n$ be a holomorphic self mapping of $\ID$, $\lambda\in[0,1)$.
Then inequality $(\ref{P5eq17})$ is true for
$$
|z|=r\leq
\begin{cases}
1/(1+2\lambda)\,\,\mbox{if}\,\, \lambda\geq 1/2,\\
\sqrt{(1/2)(1-\lambda)}\,\,\mbox{if}\,\,\lambda\leq 1/2.
\end{cases}
$$
When $\lambda\geq 1/2$, the radius $1/(1+2\lambda)$ is the best possible.
\end{Thm}
Clearly this theorem is an extension of Bohr's theorem (Theorem A) for functions with fixed initial coefficient. It is worth
mentioning that the optimal Bohr radius in Theorem~B for $\lambda<1/2$ is still unknown.
Since in Theorem~A we can consider $a_0\geq 0$ without loss of generality, varying $\lambda$ from $0$ to $1$ it is
immediately seen from Theorem~B that $(\ref{P5eq17})$ is true for all holomorphic self maps of $\ID$ whenever $r\leq 1/3$.

On the other hand, in \cite[Theorem 2.1]{Paul} Paulsen and Singh have shown that Theorem~A admits generalization for
operator valued holomorphic functions, i.e. holomorphic functions $f$ from $\ID$ to $\IB(\IH)$ where $\IB(\IH)$ is the set of bounded
linear operators on a complex Hilbert space $\IH$. To be more specific, \cite[Theorem 2.1]{Paul} was proved under the assumption
$\mbox{Re}(f(z))\leq I$, which is weaker than the natural generalization $\|f(z)\|<1$ of
the condition $f(\ID)\subset\ID$ from scalar valued case. Here $\|.\|$ denotes the operator norm in $\IB(\IH)$.
Further extensions of \cite[Theorem 2.1]{Paul} and other results of similar flavor can be found in
\cite{Pop, Pop1, Pop2}. In particular, another generalization of
$(\ref{P5eq17})$ for operator valued functions in single complex variable is available in \cite[Cor. 2.11]{Pop1} under the
assumption $f(0)=a_0I, a_0\in\IC$. 

At this point, we note that
if $f(z)=\sum_{n=0}^\infty A_nz^n$ is an operator valued holomorphic function in $\ID$, then
an operator inequality of the form $\sum_{n=0}^\infty |A_n|r^n\leq I$ also provides a perfect
analogue of the classical Bohr inequality,
$|.|$ being the absolute value of any operator in $\IB(\IH)$. Yet it appears to be less
investigated compared to the inequalities associated with norms.
Moreover, to the best of our knowledge, till date no attempt has been made to extend Theorem~B for operator valued functions.
Motivated by these facts, we prove an analogue of Theorem~B for operator valued holomorphic
functions in the following form, under more general assumptions than $A_0=a_0I, a_0\in\IC$.
\bthm\label{P5thm1}
Let $f:\ID\to\B(\IH)$ be holomorphic with an expansion $f(z)=\sum_{n=0}^\infty A_nz^n$, $A_n\in\B(\IH)$ for all
$n\in\IN\cup\{0\}$, where $A_0$ is normal and $A_0A_n=A_nA_0$ for all $n\geq 1$. If $\|f(z)\|<1$,
then
$\sum_{n=0}^\infty|A_n|r^n\leq I$ for
$$
rI\leq
\begin{cases}
(I+2|A_0|)^{-1}\,\mbox{when}\, |A_0|\geq (1/2)I,\\
((1/2)(I-|A_0|))^{1/2}\, otherwise,\\
\end{cases}
$$
where $|z|=r$.
\ethm
\br
The original proof of Theorem~B (cf. \cite{Bo} or \cite[Chapter 8]{Gar}) seems difficult to be imitated for operator valued
functions, and possibly requires stronger assumptions than that of Theorem~\ref{P5thm1}. Our proof of Theorem~\ref{P5thm1}
therefore uses operator theoretic methods combined with a few function theoretic techniques from \cite{Gol, Lit, Rogo}.
\er

Unlike the scalar valued case, an operator valued analogue of classical Bohr inequality does not follow directly from Theorem \ref{P5thm1}.
However, making use of some parts of the proof of Theorem \ref{P5thm1}, we obtain the following generalization of Theorem~A.
\bcor\label{P5cor1}
Under the hypotheses of Theorem \ref{P5thm1},
$\sum_{n=0}^\infty|A_n|r^n\leq I$
for $|z|=r\leq 1/3$.
\ecor
Another important result in the realm of Bohr radius problem is \cite[Teorema B]{Bo}
(see \cite[Theorem 8.6.15]{Gar} for the current version), which provides sharp estimate on the majorant
series $\sum_{n=0}^\infty|a_n|r^n$ of a holomorphic self mapping $f$ of $\ID$ for $|z|=r\in[1/3, 1/\sqrt{2}]$.
A part of the research on Bohr phenomenon in abstract frameworks is inspired by this result (see for example \cite{Def3}).
That said, we are unaware of any existing result in operator valued setup that directly reduces
to \cite[Theorem 8.6.15]{Gar} when restricted to scalar valued case.
In this article, the following
analogue of the aforesaid result for an operator valued holomorphic function is obtained.
\bcor\label{P5cor2}
Let $f:\ID\to\IB(\IH)$ be holomorphic with an expansion
$f(z)=a_0I+\sum_{n=1}^\infty A_nz^n$ for some $a_0\in\IC$. If $\|f(z)\|<1$ then
\be\label{P5eq15}
|a_0|I+\sum_{n=1}^\infty|A_n|r^n\leq
\left(\frac{3-\sqrt{8(1-r^2)}}{r}\right)I
\ee
for $1/3\leq |z|=r\leq 1/\sqrt{2}$.
\ecor
\br
Let $f(z)=\sum_{n=0}^\infty A_nz^n$ be a $\IB(\IH)$
valued holomorphic function on $\ID$ satisfying $\|f(z)\|<1$.
From \cite{BB2} it is already known that
$\sum_{n=0}^\infty|A_n|r^n\leq (1/\sqrt{1-r^2})I$ for any $|z|=r<1$. Since
$$
1\leq\frac{3-\sqrt{8(1-r^2)}}{r}
\leq \frac{1}{\sqrt{1-r^2}}
$$
for any $1/3\leq r\leq 1/\sqrt{2}$, the above two corollaries improve that result from \cite{BB2} in
certain ranges of $r$, but under additional hypotheses.
\er

In view of the above results, it is natural to ask if the
functions satisfying\,\, $\mbox{Re}(f(z))\leq I$ for each $z\in\ID$ exhibit Bohr phenomenon
involving the absolute values of operators. It turns out that if we consider $A_0=a_0I, a_0\geq 0$, then \cite[Theorem 2.1(6)]{Paul}
readily gives one such Bohr inequality. In the following theorem we prove another result of this kind.
\bthm\label{P5thm2}
Let $f:\ID\to\B(\IH)$ be holomorphic with an expansion $f(z)=\sum_{n=0}^\infty A_nz^n$, $A_n\in\B(\IH)$ for all
$n\in\IN\cup\{0\}$, and $\emph{Re}(f(z))\leq I$ for all $z\in\ID$.
In addition, we assume that $f(z)$ is normal for each $z\in\ID$, $A_0A_n=A_nA_0$ for all $n\geq 1, A_0\geq 0$ and $\|A_0\|<1$.
Then $\sum_{n=0}^\infty|A_n|r^n\leq I$ for all $|z|=r\leq 1/3$.
\ethm
\section{Proofs of the main results}
We list down the following known results which we will use frequently in our proofs.

\vspace{3pt}
\noindent\textit{\bf{Maximum Principle [MP]}}(cf. \cite[p. 3]{Goh}): Let $D\subset\IC$ be a bounded open set. Also let $E$ be a complex Banach space and
let $f: \ov{D}\to E$ be a continuous function which is holomorphic in $D$. Denote by $\partial D$ the boundary of $D$. Then
$$
\max_{z\in\ov{D}}\|f(z)\|=\max_{z\in\partial D}\|f(z)\|.
$$

\vspace{3pt}
\noindent\textit{\bf{Fuglede-Putnam Theorem [FPT]}}(see f.i. \cite[p. 278]{Con}): If $N$ and $M$ are normal operators on the complex Hilbert spaces
$\IH$ and $\K$ and $B:\K\to\IH$ is an operator such that $NB=BM$, then $N^*B=BM^*$.

\vspace{3pt}
\noindent\textit{\bf{Properties of Positive Operators [PPO]}}(cf. \cite[pp. 261--262]{Fri}): Let $A, B\in\B(\IH)$ be such that $A\geq 0$, $B\geq 0$ and
$AB=BA$. Then
\bee
\item[(i)]~
$A^{1/2}B^{1/2}=B^{1/2}A^{1/2}.$
\item[(ii)]~
$AB\geq 0$ and $(AB)^{1/2}=A^{1/2}B^{1/2}$.
\item[(iii)]~
If $A$ is invertible then $A^{-1}\geq 0$. Also $A^{1/2}$ is invertible and
$(A^{1/2})^{-1}=(A^{-1})^{1/2}$.
\item[(iv)]~
If $A\geq B$, $A$ and $B$ are invertible then $B^{-1}\geq A^{-1}$.
\item[(v)]~
If $A\geq B$ then $A^{1/2}\geq B^{1/2}$.
\eee
We are now ready to prove the first theorem of this article.
\bpf[Proof of Theorem \ref{P5thm1}]
Since $A_0A_n=A_nA_0$ for $n\geq 1$, $A_0f(z)=f(z)A_0$ for each $z\in\ID$. As $A_0$ is normal, using [FPT] repeatedly we obtain
\be\label{P5eq3}
A_0^*f(z)=f(z)A_0^*\implies f(z)^*A_0=A_0f(z)^*\implies f(z)^*A_0^*=A_0^*f(z)^*,
\ee
and similarly,
\be\label{P5eq4}
A_0^*A_n=A_nA_0^*\implies A_n^*A_0=A_0A_n^*\implies A_n^*A_0^*=A_0^*A_n^*.
\ee
From (\ref{P5eq3}) it is easy to see that $|f(z)|^2|A_0|^2=|A_0|^2|f(z)|^2$, and hence
$$
(I-|f(z)|^2)(I-|A_0|^2)=(I-|A_0|^2)(I-|f(z)|^2)
$$
for any $z\in\ID$. Observing
that both $I-|f(z)|^2$ and $I-|A_0|^2$ are positive operators,
[PPO(ii)] yields $(I-|f(z)|^2)(I-|A_0|^2)\geq 0$, which is same as saying
$$
I-A_0^*f(z)-f(z)^*A_0+f(z)^*f(z)A_0^*A_0
\geq f(z)^*f(z)-A_0^*f(z)-f(z)^*A_0+A_0^*A_0.
$$
As $A_0$ commutes with both $f(z)$ and $A_0^*$, from the above
inequality we get
$$
(I-A_0^*f(z))^*(I-A_0^*f(z))\geq (f(z)-A_0)^*(f(z)-A_0),
$$
or, in other words
\be\label{P5eq5}
\|(I-A_0^*f(z))x\|\geq\|(f(z)-A_0)x\|
\ee
for any $x\in\IH$. Since $\|A_0^*f(z)\|\leq \|A_0\|\|f(z)\|<1$,
$I-A_0^*f(z)$ is invertible, and hence by virtue of (\ref{P5eq5}),
\be\label{P5eq13}
\phi(z):=(f(z)-A_0)(I-A_0^*f(z))^{-1}
\ee
is a holomorphic function from $\ID$ to $\IB(\IH)$ with
$\|\phi(z)\|\leq 1$ and $\phi(0)=0$. Again, due to the fact that
$\|A_0^*\phi(z)\|<1$, $(I+A_0^*\phi(z))$ is invertible. Also using (\ref{P5eq3}) suitably,
we get that $\phi(z)$ commutes with both $A_0$ and $A_0^*$,
and hence from (\ref{P5eq13})
$$
f(z)=(A_0+\phi(z))(I+A_0^*\phi(z))^{-1},
$$
or equivalently
\be\label{P5eq6}
\sum_{n=0}^\infty A_nz^n=
A_0+(I-|A_0|^2)\sum_{n=1}^\infty{(-A_0^*)}^{n-1}{\phi}^n(z).
\ee
Therefore for any fixed $k\geq 1$,
\be\label{P5eq8}
A_0+(I-|A_0|^2)\sum_{n=1}^k{(-A_0^*)}^{n-1}{\phi}^n(z)
=\sum_{n=0}^kA_nz^n+\sum_{n=k+1}^\infty C_nz^n
\ee
for some bounded linear operator $C_n$'s. Let us define $g$ by
$$
g(z)=A_0+(I-|A_0|^2)\sum_{n=1}^k{(-A_0^*)}^{n-1}z^n,\, z\in\ID.
$$
Further, we write $g(z)=\sum_{n=0}^k G_nz^n$,
where $G_n$ is the $n^{th}$ $\IB(\IH)$ valued coefficient of $g$,
i.e. $G_0=A_0$, $G_n=(I-|A_0|^2){(-A_0^*)}^{n-1}$ for $1\leq n\leq k$
and $G_n=0$ for $n\geq k+1$.
Now we define $g_i$, $i\in\IN\cup\{0\}$ by
$$
g_i(z)=\sum_{n=0}^k G_{n+i}z^n,\,z\in\ID.
$$
Clearly $g_0(z)=g(z)$ and $g_i\equiv 0$ for any $i\geq k+1$.
Also let
$$
g_i[\phi(z)]:=\sum_{n=0}^k G_{n+i}\phi^n(z)=\sum_{n=0}^\infty G_n(i)z^n,\,G_n(i)\in\IB(\IH)
$$
for any $i\in\IN\cup\{0\}$.
Again it is evident that for any $n\geq0$, $G_n(i)=0$ if $i\geq k+1$.
Since $\phi$ commutes with any $G_n$, $0\leq n\leq k$, we have
$g_i[\phi(z)]=g_i(0)+\phi(z)g_{i+1}[\phi(z)]$ for any $0\leq i\leq k$.
As $\phi(0)=0$, for any $i$ we write
$$
\phi(z)g_{i}[\phi(z)]=\sum_{n=1}^\infty P_n(i)z^n,
$$
where $P_n(i)\in\IB(\IH)$.
Now using the fact $\|\phi(z)\|\leq 1$, it is readily seen that for any $x\in\IH$
\begin{equation*}
\begin{split}
\MoveEqLeft[-0.5]
\|g_i[\phi(z)]x\|^2=\|(g_i(0)+\phi(z)g_{i+1}[\phi(z)])x\|^2\\
&=\|g_i(0)x\|^2+\|\phi(z)g_{i+1}[\phi(z)]x\|^2+
2\mbox{Re}\langle\phi(z)g_{i+1}[\phi(z)]x, g_i(0)x\rangle\\
&\leq\|G_ix\|^2+\left\|\left(\sum_{n=0}^\infty G_n(i+1)z^n\right)x\right\|^2\hspace{-5pt}+
2\mbox{Re}\left\langle\left(\sum_{n=1}^\infty P_n(i+1)z^n\right)x, g_i(0)x\right\rangle.
\end{split}
\end{equation*}
Putting $z=re^{i\theta}$ in the above inequality and then integrating both sides of this inequality
over $\theta$ from $0$ to $2\pi$, we obtain
\be\label{P5eq7}
\sum_{n=0}^\infty|G_n(i)|^2r^{2n}
\leq|G_i|^2+\sum_{n=0}^\infty|G_n(i+1)|^2r^{2n},
\ee
because
\begin{equation*}
\begin{split}
\MoveEqLeft[-2]
\int_{\theta=0}^{2\pi}2\mbox{Re}\left\langle\left(\sum_{n=1}^\infty P_n(i+1)r^ne^{in\theta}\right)x, g_i(0)x\right\rangle d\theta\\
&=\sum_{n=1}^\infty r^n\left(\langle P_n(i+1)x, g_i(0)x\rangle\int_{\theta=0}^{2\pi}e^{in\theta}d\theta+
\langle g_i(0)x, P_n(i+1)x\rangle\int_{\theta=0}^{2\pi}e^{-in\theta}d\theta\right)
=0.
\end{split}
\end{equation*}
As a result, summing both sides of the inequality $(\ref{P5eq7})$ from $i=0$ to $k$ and then allowing $r\to 1-$, we get
$$
\sum_{n=0}^\infty|G_n(0)|^2\leq\sum_{i=0}^k|G_i|^2+\sum_{n=0}^\infty|G_n(k+1)|^2=\sum_{i=0}^k|G_i|^2,
$$
which is indeed same as saying (see (\ref{P5eq8}))
$$
\sum_{n=0}^k|A_n|^2+\sum_{n=k+1}^\infty|C_n|^2\leq
\sum_{i=0}^k|G_i|^2=|A_0|^2+(I-|A_0|^2)^2
\sum_{n=1}^k|A_0^*|^{2(n-1)}.
$$
As $A_0$ is normal, $|A_0|=|A_0^*|$, and therefore the above inequality implies
\be\label{P5eq9}
\sum_{n=1}^k|A_n|^2\leq(I-|A_0|^2)^2\sum_{n=1}^k|A_0|^{2(n-1)}.
\ee
Now from (\ref{P5eq4}) one can prove that, $|A_0|^2|A_n|^2=|A_n|^2|A_0|^2$ for any $n\geq 1$, and
hence by [PPO(i)], $|A_0||A_n|=|A_n||A_0|$, which therefore implies that
$|A_0|^p|A_n|^q=|A_n|^q|A_0|^p$ for any $p,q\in\IN\cup\{0\}$. Using this fact, a little computation reveals that
\begin{align*}
\left(\sum_{n=1}^k|A_n|^2\right)
\hspace{-3pt}\left(\sum_{n=1}^k|A_0|^{2n}\right)-\left(\sum_{n=1}^k|A_n||A_0|^n\right)^2
\hspace{-8pt}=\sum_{i=1}^{k-1}\sum_{j=i+1}^k
\left(|A_i||A_0|^j-|A_j||A_0|^i\right)^2\geq 0,
\end{align*}
as for any $i, j$, $|A_i||A_0|^j-|A_j||A_0|^i$ is self adjoint.
Now combining (\ref{P5eq9}) and the previous inequality, we have,
by [PPO(ii)]
$$
\left(\sum_{n=1}^k|A_n||A_0|^n\right)^2
\leq|A_0|^2(I-|A_0|^2)^2\left(\sum_{n=1}^k|A_0|^{2(n-1)}\right)^2,
$$
or, after taking the square root on both sides,
\be\label{P5eq10}
\sum_{n=1}^k|A_n||A_0|^n\leq
|A_0|(I-|A_0|^2)\sum_{n=1}^k|A_0|^{2(n-1)}.
\ee
Assuming $|A_0|$ invertible, we
observe that for any $n\in\IN$
and for any fixed $r\in(0,1)$,
$$
r^n|A_0|^{-n}\geq r^{n+1}|A_0|^{-(n+1)}\iff rI\leq|A_0|.
$$
Therefore,
if $rI\leq|A_0|$, then doing a summation by parts
and applying (\ref{P5eq10}) through
[PPO(ii)], we get
\begin{equation*}
\begin{split}
\MoveEqLeft[-1]
\sum_{n=1}^k|A_n|r^n=\sum_{n=1}^k(|A_n||A_0|^n)(r^n|A_0|^{-n})\\
&=\sum_{n=1}^{k-1}\left(r^n|A_0|^{-n}-r^{n+1}|A_0|^{-(n+1)}\right)
\left(\sum_{t=1}^n|A_t||A_0|^t\right)
+r^k|A_0|^{-k}\left(\sum_{n=1}^k|A_n||A_0|^n\right)\\
&\leq\sum_{n=1}^{k-1}\left(r^n|A_0|^{-n}-r^{n+1}|A_0|^{-(n+1)}
\right)
\left(|A_0|(I-|A_0|^2)\sum_{t=1}^n|A_0|^{2(t-1)}
\right)\\
&+r^k|A_0|^{-k}\left(|A_0|(I-|A_0|^2)\sum_{n=1}^k|A_0|^{2(n-1)}
\right)=|A_0|(I-|A_0|^2)\sum_{n=1}^k|A_0|^{n-2}r^n,
\end{split}
\end{equation*}
and consequently,
\be\label{P5eq11}
\sum_{n=1}^\infty|A_n|r^n
\leq(I-|A_0|^2)\sum_{n=1}^\infty|A_0|^{n-1}r^n
=r(I-|A_0|^2)(I-r|A_0|)^{-1}
\ee
for $rI\leq|A_0|$.
Further, using [PPO(ii), (iii)] appropriately, from inequality $(\ref{P5eq11})$ we find that
$\sum_{n=0}^\infty|A_n|r^n\leq I$ is satisfied if
$|A_0|$ is invertible and the two conditions
$$
rI\leq |A_0|\,,\,r(I+2|A_0|)\leq I
$$
hold together. Now
when $|A_0|$ is invertible,
$\||A_0|x\|\geq\delta\|x\|$ for some $\delta>0$ and for all $x\in\IH$,
which is same as saying $|A_0|\geq\delta I$, or equivalently $\langle |A_0|x, x\rangle\geq\delta\langle x,x\rangle$
for any $x\in\IH$. Hence
$$
\|(I+2|A_0|)x\|\|x\|\geq\langle (I+2|A_0|)x, x\rangle\geq\|x\|^2+2\delta\|x\|^2,
$$
i.e. $\|(I+2|A_0|)x\|\geq(1+2\delta)\|x\|$ for any $x\in\IH$. Since
$I+2|A_0|$ is positive, the previous inequality implies that $I+2|A_0|$ is invertible.
Thus by [PPO(iii)], $(I+2|A_0|)^{-1}\geq 0$. Therefore
$$
r(I+2|A_0|)\leq I \iff rI\leq (I+2|A_0|)^{-1}
$$
(cf. [PPO(ii)]). Now suitably using [PPO(ii), (iii)]
and after a little calculation, we get
$$
(I+2|A_0|)^{-1}\leq|A_0|\iff |A_0|\geq (1/2)I,
$$
i.e. $r(I+2|A_0|)\leq I$ and $rI\leq|A_0|$ are simultaneously satisfied if $|A_0|\geq (1/2)I$, which also implies the
invertibility of $|A_0|$. Therefore if $|A_0|\geq(1/2)I$, then
$$
rI\leq(I+2|A_0|)^{-1}\implies\sum_{n=0}^\infty|A_n|r^n\leq I.
$$

On the other hand, $\|f(z)\|<1$ implies
$\|f(z)x\|^2\leq\|x\|^2$ for all $x\in\IH$. Therefore setting $z=re^{i\theta}$ and integrating
both sides of this inequality over $\theta$ from $0$ to $2\pi$, and then letting $r\to1-$,
we have $\sum_{n=1}^\infty|A_n|^2\leq I-|A_0|^2$. Now
\cite[Lemma1]{BB2} gives
\be\label{P5eq12}
\sum_{n=1}^\infty|A_n|r^n\leq(I-|A_0|^2)^{1/2}
\frac{r}{\sqrt{1-r^2}}
\ee
for any $r\in[0,1)$. Therefore, $\sum_{n=0}^\infty|A_n|r^n\leq I$ whenever
$$
(I-|A_0|^2)^{1/2}\frac{r}{\sqrt{1-r^2}}\leq I-|A_0|.
$$
Since $(I-|A_0|^2)^{1/2}$ and $(I-|A_0|)$ both are positive and commute with each other, the above
inequality is equivalent to
$$
(I-|A_0|^2)\frac{r^2}{1-r^2}\leq (I-|A_0|)^2.
$$
A use of [PPO(ii), (iii)] reveals that
this inequality can be reduced to
$$
(I+|A_0|)r^2\leq(I-|A_0|)(1-r^2),
$$
i.e. $rI\leq((1/2)(I-|A_0|))^{1/2}$. Hence
$$
rI\leq((1/2)(I-|A_0|))^{1/2}\implies\sum_{n=0}^\infty|A_n|r^n\leq I
$$
without any extra assumption on $|A_0|$, and thus it remains true if $|A_0|-(1/2)I$ is not a
positive operator.
This completes our proof.
\epf

Using some facts from the above proof, we now establish the two corollaries.
\bpf[Proof of Corollary \ref{P5cor1}]
From the proof of Theorem \ref{P5thm1},
we observe that $\tilde{\phi}(z):=\phi(z)/z$ is an operator valued holomorphic
function in $\ID$, where $\phi$ is given by $(\ref{P5eq13})$. By [MP], it is immediate that
for any $r\in(0,1)$, $\|\tilde{\phi}(z)\|\leq 1/r$ for $|z|\leq r$.
Letting $r\to 1-$, we have $\|\tilde{\phi}(z)\|\leq 1$ for all $z\in\ID$. In particular for $z=0$,
we get $\|A_1(I-|A_0|^2)^{-1}\|\leq1$, i.e. $\|A_1x\|\leq\|(I-|A_0|^2)x\|$ for any $x\in\IH$. Thus we get the following:
\be\label{P5eq14}
|A_1|\leq I-|A_0|^2\leq 2(I-|A_0|).
\ee
Now for any given $n>1$, let $w$ be a primitive $n^{th}$ root
of unity. Then we observe that the function
$$
F_1(z):=\frac{1}{n}\sum_{k=0}^{n-1} f(w^{k}z)=A_0+A_nz^n+A_{2n}
z^{2n}+\cdots, z\in\ID
$$ satisfies all the assumptions we began with, and so does the function $F(z)=A_0+A_nz+A_{2n}z^{2}+\cdots, z\in\ID$.
Therefore the inequalities in (\ref{P5eq14}) remain true if we replace $A_1$ by
$A_n$. Rest of the proof follows from straightforward calculations.
\epf
\bpf[Proof of Corollary \ref{P5cor2}]
This proof relies on a computational trick from the proof of \cite[Theorem 2.1(4)]{Paul}.
It is evident that for any fixed $r\in[0,1)$, either $|a_0|I\geq rI$ or $|a_0|I\leq rI$. Hence, looking at
(\ref{P5eq11}) and (\ref{P5eq12}), it becomes clear that for any $r\in(1/3, 1/\sqrt{2}]$
\be\label{P5eq16}
|a_0|I+\sum_{n=1}^\infty |A_n|r^n\leq
\max\{\max_{x\in[0,r]}\chi(x),\max_{x\in[r,1)}\xi(x)\}I,
\ee
where
$$
\chi(x):=x+r\sqrt{1-x^2}/\sqrt{1-r^2}\,\,
\mbox{and}\,\, \xi(x):=x+r(1-x^2)/(1-rx)
$$
for $x\in[0,1)$.
Now
$\xi(x)$ is strictly increasing for $x<x_0$ and stricly
decreasing for $x>x_0$, where $x_0=(1/r)(1-\sqrt{(1/2)(1-r^2)})$.
Note that $r\leq x_0< 1$ for $r\in(1/3,1/\sqrt{2}]$. Therefore
$$
\max_{x\in[r,1)}\xi(x)=\xi(x_0)=\left(3-\sqrt{8(1-r^2)}\right)\Big/r.
$$
On the other hand, $\chi(x)$ is strictly increasing for $x<\sqrt{1-r^2}$ and strictly
decreasing for $x>\sqrt{1-r^2}$. For $r\leq1/\sqrt{2}$, $r\leq \sqrt{1-r^2}$; and hence
$$
\max_{x\in[0,r]}\chi(x)=\chi(r)=2r.
$$
Observing that
$$
r\in(0,1)\implies 2r\leq\left(3-\sqrt{8(1-r^2)}\right)\Big/r,
$$
and that $\left(3-\sqrt{8(1-r^2)}\right)\Big/r=1$ at
$r=1/3$, our proof is complete from (\ref{P5eq16}).
\epf

We now prove the final result of this article.
\bpf[Proof of Theorem \ref{P5thm2}]
It is immediately seen that for any $x\in\IH$,
\begin{equation*}
\begin{split}
\MoveEqLeft[-0.5]
\|(2(I-A_0)-(f(z)-A_0))x\|^2
=4\|(I-A_0)x\|^2\\
&-4\langle \mbox{Re}((I-A_0)(f(z)-A_0))x, x\rangle+ \|(f(z)-A_0)x\|^2.
\end{split}
\end{equation*}
Now since $A_0A_n=A_nA_0$ for all $n\in\IN$ and $f(0)=A_0$ is normal,
by the same argument as in the proof of Theorem \ref{P5thm1},
we have $A_0f(z)=f(z)A_0$ and $A_0f(z)^*=f(z)^*A_0$ for each $z\in\ID$ (see (\ref{P5eq3})). Using this fact,
and that $\mbox{Re}(f(z))\leq I, A_0=A_0^*$ we have
\begin{equation*}
\begin{split}
\MoveEqLeft[-0.5]
\mbox{Re}((I-A_0)(f(z)-A_0))\\
&=(I-A_0)\mbox{Re}(f(z)-A_0)\\
&=(I-A_0)(\mbox{Re}(f(z)-A_0)-(I-A_0))+(I-A_0)^2\leq (I-A_0)^2
\end{split}
\end{equation*}
(see [PPO(ii)]), which implies $\langle\mbox{Re}((I-A_0)(f(z)-A_0))x, x\rangle\leq\|(I-A_0)x\|^2$
for any $x\in\IH$. Therefore for each $x\in\IH$
\be\label{P5eq1}
\|(2(I-A_0)-(f(z)-A_0))x\|\geq\|(f(z)-A_0)x\|,
\ee
which further gives
\begin{equation*}
\begin{split}
\MoveEqLeft[6]
2\|(I-A_0)x\|\leq\|(2(I-A_0)-(f(z)-A_0))x\|+\|(f(z)-A_0)x\|\\
&\leq 2\|(2(I-A_0)-(f(z)-A_0))x\|.
\end{split}
\end{equation*}
Since $\|A_0\|<1$, $I-A_0$ is invertible and hence
bounded below. Therefore $2(I-A_0)-(f(z)-A_0)$ is also bounded
below by the above inequality.
Observing that $f(z)-A_0$ and $I-A_0$ both are normal operators
that commute with each other, $2(I-A_0)-(f(z)-A_0)$ is again normal, and therefore invertible.
As a result, from (\ref{P5eq1}) we find that
$$
g_1(z):=(f(z)-A_0)(2(I-A_0)-(f(z)-A_0))^{-1},\,z\in\ID
$$
is an operator valued holomorphic function with $\|g_1(z)\|\leq 1$ for each $z\in\ID$.
Again $\phi_1(z):=g_1(z)/z, z\in\ID$ is
an operator valued holomorphic function, and hence for any $r\in(0,1)$ we have,
by [MP], that $\|\phi_1(z)\|\leq 1/r$ for all $|z|\leq r$. Letting $r\to1-$ we conclude that
$\|\phi_1(z)\|\leq 1$ for all $z\in\ID$.
Now
$$
(f(z)-A_0)=z\phi_1(z)(2(I-A_0)-(f(z)-A_0)),
$$
and hence
$$
\|(f(z)-A_0)x\|^2\leq |z|^2\|(2(I-A_0)-(f(z)-A_0))x\|^2.
$$
Putting $z=re^{i\theta}$ and inserting the expansion of $f(z)$
in the above inequality, we get
$$
\left\|\left(\sum_{n=1}^\infty A_nr^ne^{in\theta}\right)x\right\|^2
\leq\left\|\left(2(I-A_0)re^{i\theta}-\sum_{n=1}^\infty A_nr^{n+1}e^{i(n+1)\theta}\right)x\right\|^2.
$$
Now integrating both sides of this inequality over $\theta$ from
$0$ to $2\pi$, we obtain
$$
\sum_{n=1}^\infty|A_n|^2r^{2n}\leq 4(I-A_0)^2r^2+\sum_{n=1}^\infty|A_n|^2r^{2n+2}
$$
for any $r\in[0,1)$. Replacing $r^2$ by $r$, a little calculation yields that
\be\label{P5eq2}
\sum_{n=1}^\infty|A_n|^2r^{n}\leq 4(I-A_0)^2\frac{r}{1-r}.
\ee
Setting $H_n=A_nr^{n/2}, k=1$ and replacing $r$ by $\sqrt{r}$ in \cite[Lemma 1]{BB2},
a use of the inequality (\ref{P5eq2}) now reveals that
$$
\sum_{n=1}^\infty|A_n|r^n\leq 2(I-A_0)\frac{r}{1-r}
$$
for each $r\in[0,1)$. The desired result follows immediately from the above inequality.
\epf
\section{Concluding remarks}
It is worth noting that both Theorem~A and Theorem~B can be proved under the weaker
hypothesis $|f(z)|\leq 1$, as in this situation $|f(z)|=1$ for any $z\in\ID$
implies, by the Maximum modulus principle of complex valued functions, that $f(z)=1$ for all $z\in\ID$.
An exact analogue of this principle is not available for operator valued functions in general.
In view of that, it is natural to ask if
the condition $\|f(z)\|<1$ in the statements of Theorem \ref{P5thm1} and Corollaries \ref{P5cor1} and \ref{P5cor2}
could be replaced by $\|f(z)\|\leq 1$. This appears to be an interesting problem for future research.

\end{document}